\documentclass[11pt,reqno]{amsart}

\usepackage{latexsym}
\usepackage{amsfonts}
\usepackage{amsmath}
\usepackage{url}
\usepackage{graphicx}
\usepackage{color}
\usepackage{bbm}
\definecolor{myblue}{rgb}{0.09,0.32,0.44} %22-84-113
\usepackage[pdfborder={0 0 0},pdfborderstyle={},colorlinks=true,linkcolor=myblue,citecolor=myblue,urlcolor=blue]{hyperref}

\usepackage{epsfig,psfrag,verbatim,amsfonts}

\usepackage{amssymb,amsmath}

\def\bl{\begin{lemma}}
\def\el{\end{lemma}}
\def\bth{\begin{theorem}}
\def\eth{\end{theorem}}
\def\bc{\begin{corollary}}
\def\ec{\end{corollary}}
\def\bcj{\begin{conjecture}}
\def\ecj{\end{conjecture}}
\def\bpr{\begin{proposition}}
\def\epr{\end{proposition}}
\def\bde{\begin{definition}}
\def\ede{\end{definition}}

\newcommand{\be}{\begin{eqnarray}}
\newcommand{\ee}{\end{eqnarray}}

\renewcommand{\and}{\hbox{ {\rm and} }}

\newcommand{\LL}{L}

\newtheorem{theorem}{Theorem}[section]
\newtheorem{definition}{Definition}[section]
\newtheorem{lemma}[theorem]{Lemma}

\newtheorem{corollary}[theorem]{Corollary}
\newtheorem{proposition}[theorem]{Proposition}
\newtheorem{conjecture}[theorem]{Conjecture}

\newtheorem*{theorem*}{Theorem}

\theoremstyle{definition}
\newtheorem*{definition*}{Definition}
\numberwithin{equation}{section}

\begin{document}
\title{continuous vs discrete spins in the hyperbolic plane}
\author{Itai Benjamini \and Gady Kozma }

\date{March 2016}

\begin{abstract}
We study the $O(n)$ model on graphs quasi-isometric to the hyperbolic plane, with free boundary conditions. We observe that
%While it is known that the (free) Ising model, $O(1)$, admits a second order phase transition for lattices in the hyperbolic plane.
%We observe that when the spins take continuous values, as in the $O(n)$ models for $n > 1$,
the pair correlation decays exponentially with distance, for all temperatures, if and only if $n>1$.
\end{abstract}

\maketitle

%\section{Introduction}

We wish to report here on a curious contrast between the behaviour of models with discrete and continuous symmetry on graphs quasi-isometric to the hyperbolic plane.
For concreteness we will study the $O(n)$ model, with $n=1$, or Ising, being our model for a model with discrete spins and $n\ge 2$ our model for continuous spins and symmetry group. Let us start with the relevant definitions.
\begin{definition*}Two metric spaces $X$ and $Y$ are called {\em quasi-isometric} if there is a map $\phi:X\to Y$ with the properties
\begin{gather*}
\forall x,y\in X\qquad cd(x,y)-C\le d(\phi(x),\phi(y))\le Cd(x,y)+C\\
\forall y\in Y\exists x\in X \textrm{ such that }d(\phi(x),y)\le C
\end{gather*}
where $c$ and $C$ denote constants which are independent of $x$ and $y$.
\end{definition*}
We say that a graph is quasi-isometric to the hyperbolic plane if, when you equip it with the graph metric (i.e.\ the distance between any two vertices is the length of the shortest path between them), it is quasi-isometric to the hyperbolic plane with its standard metric.

An easy example of a graph quasi-isometric to the hyperbolic plane is the seven-regular planar triangulation: the (unique) graph which is the 1-skeleton of an infinite triangulation of a planar disk which satisfies that the degree of every vertex is 7 (the number 7 may be replaced with any $n>6$ and the graph would still be quasi-isometric to the hyperbolic plane). This graph may be constructed completely combinatorially, if desired. Another example can be had by taking a binary tree and connecting each generation ``horizontally''. See \cite[\S 14]{CFLP}. Another family of examples can be achieved from cocompact lattices of $\textrm{PSL}_2(\mathbb{R})$ --- the action of $\textrm{PSL}_2(\mathbb{R})$ on the hyperbolic plane by M\"obius transformations gives, for each such lattice, a uniformly discrete subset of the hyperbolic plane, over which one may impose a graph structure in various ways, e.g.\ by taking the Voronoi triangulation of it. See \cite{CFLP} for an introduction to hyperbolic geometry and \cite{K} for lattices in $\textrm{PSL}_2(\mathbb{R})$.

Let us also recall the definition of the $O(n)$ model. We endow every vertex $v$ of our graph  with a spin $s_v$, which are unit vectors in the $(n-1)$-sphere $\mathbb{S}^{n-1}$. The Hamiltonian is given by:
$$
H = -\sum_{\{v,u\}\in E}\langle s_v,s_u\rangle .
$$
We take an exhaustion of our graph, say the balls in the graph distance $B_r$, and define the pair correlation as
\[
\lim_{r\to\infty}\int\langle s_u,s_v\rangle \exp(-\beta H_r(s))\,ds\bigg/\int \exp(-\beta H_r(s))\,ds
\]
where the integrals are taken over $(\mathbb{S}^{n-1})^{|B_r|}$ and $H_r$ is the Hamiltonian with the sum restricted to edges inside $B_r$. One may also wire all the vertices of the boundary of $B_r$ to get the Hamiltonian with wired boundary conditions, and the corresponding pair correlation.

As stated, our purpose in this note is to compare the pair correlation for $n=1$ and $n\ge 2$ on graphs quasi-isometric to the hyperbolic plane.
Let us start with the case of {\em discrete spins}, specifically Ising (or $n=1$), which is well-known. In \cite{BS1} it was shown, using a Peierls-type argument  that the threshold for uniqueness of the infinite cluster, $p_u$, is strictly smaller than $1$ for Bernoulli percolation
on any transitive planar nonamenable graph with one end. In particular this applies to graphs quasi-isometric to the hyperbolic plane, as both amenability and the number of ends are invariant to quasi-isometry. For any graph this implies that the Ising model has a phase transition for the pair correlation, see \cite[\S 4]{ACCN88} (the argument of \cite{ACCN88} is sufficiently simple to sketch here: since the functions $\mathbbm{1}\{u\leftrightarrow v\}$ and $2^{\#\{\textrm{clusters}\}+\#\{\textrm{open edges}\}}$ are both increasing, the FKG inequality allows to control $u\leftrightarrow v$ in the FK representation at $2p/(1+p)$ which is the same as the Ising model at $\beta=|\log(2p/(1+p))|$).
%, and there is no unique infinite cluster at $p_u$.
%The argument of \cite{BS1} is a Peierls-type argument and therefore applies also to the Ising model on such graphs.
We get that the Ising model has an ordered phase, where pair correlations do not decay. This does not depend on boundary conditions (free or wired, for example). See also \cite{W} for the Ising model on hyperbolic graphs.
%By the FK presentation ***I didn't understand that, Itai, could you explain the relation to FK? I though that the relation between percolation and Ising is not so simple, that it's in fact a result of Schonmann, was I wrong?*** it follows that the (free) Ising model admits a second order phase transition on cocompact hyperbolic lattices.
%\medskip

To study {\em continuous spins}, recall that one may view a graph as an electrical network, where each edge is a $1$ Ohm conductor. The effective resistance between two vertices can then be calculated by solving Kirchoff's equations. For infinite graphs the solution is not necessarily unique, but solutions with given boundary conditions (say free or wired) are well-defined and unique. See \cite{DS84} for the necessary background.
%We show below that if the electric resistance between any pair of vertices is proportional to their graph distance,
%then pair correlation, in (free) $O(n)$ models for any $n >1$,
%decays exponentially with distance, for all temperatures. And if the resistance grows to infinity with distance, then there is no phase of magnetisation.
%\medskip
For graphs quasi-isometric to the hyperbolic plane, endowed with free boundary conditions, the electric resistance is proportional to the distance \cite{B}. This implies that for any two vertices $x$ and $y$ and for any $r$ sufficiently large one may find a function $a$ defined on $B(x,r)$, the ball of radius $r$ around $x$ in the graph distance, with the properties that
\begin{gather*}
a(y)-a(x)>c_1d(x,y)\qquad\qquad \sum_{u\sim v}|a(u)-a(v)|^2\le \tfrac12 |a(y)-a(x)|\\
\forall u\sim v \quad|a(u)-a(v)|\le \tfrac1{10}
\end{gather*}
where $c_1$ is some constant, $u\sim v$ means that $u$ and $v$ are neighbours in the graph (and are implicitly assumed to both be in $B(x,r)$), and $d(x,y)$ is the graph distance. See again \cite{DS84} for the various equivalent definitions of electrical resistance.

The existence of such an $a$ gives that the pair correlations of the $O(n)$ model, $n\ge 2$, decay exponentially, using the argument of McBryan and Spencer \cite{MS77}. The proof in that paper carries through with no change. Let us sketch the argument of \cite{MS77} as well (even though \cite{MS77} is very short and to the point). It is a version of the Mermin-Wagner argument where quadratic control of the error is achieved using a \emph{complex change of variables}. Specifically, the integral from the definition of the pair correlation is written as
\[
\int \exp\big(i(\theta(x)-\theta(y))-\beta H(\vec{\theta})\big)\qquad \theta(u)\in [0,2\pi)\forall u
\]
and then the change of variables $\theta(x)\mapsto \theta(x)+i a(x)$ is applied. The first term makes one ``earn'' $\exp(a(y)-a(x))$ while the term containing the Hamiltonian makes one lose approximately $\exp\Big(\sum_{u\sim v}|a(u)-a(v)|^2\Big)$.

We thus arrive at the following observation:

\begin{theorem*}The $O(n)$ model, $n\ge 2$ on a graph quasi-isometric to the hyperbolic plane with free boundary conditions has exponentially decaying pair correlations.\end{theorem*}

A similar result can be obtained under weaker expansion properties of the graph. In \cite{BS:har} it was shown that every transient bounded-degree planar graph admits a non-constant harmonic function with gradient in $\LL^2$.
This implies that the resistance between a pair of vertices grows to infinity with the distance \cite{BLPS, LP}. For recurrent graph this holds as well.
We conclude that the pair correlations decay for $O(n)$, $n\ge 2$ models on any bounded degree planar graph.

We remark that one cannot get such an example with trees. Take for concreteness the binary tree. Certainly, it also has resistance linear in the distance, so the $O(n)$, $n\ge 2$ has exponential decay of correlations. But the same happens at $n=1$, as for the Ising model it is the number of ends that determines the behaviour, and the tree has infinitely many ends. Thus on the tree we cannot see the change of behaviour from $n=1$ to $n \ge 2$ we see on the hyperbolic plane.

One is tempted to conjecture that for wired boundary conditions the behaviour is different and there is a phase transition, because the electric resistance is bounded.
Similarly, the electric resistance (free or wired) between a pair of vertices in lattices in the hyperbolic space (i.e.\ in dimension 3 or higher) is uniformly bounded. Thus we expect a phase transition for all $n$.

Going back to the hyperbolic plane with free boundary conditions, we do not know if there is a Kosterlitz-Thouless transition (see e.g.\ \cite{FS81,A94,F12} for a Kosterlitz-Thouless transition). Nor do not know if there is a phase transition in the total magnetisation --- standard arguments give that for temperature sufficiently high, the total magnetisation is finite. But the graph has exponential volume growth which might win over the exponential decay of the correlations at low temperatures.

Another interesting variation is to take the spins in various subsets of $\mathbb{S}^1$. For example, when the spins take the value in two intervals, is it true that one has a phase transition for the choice of the interval, but exponential decay of correlations inside the intervals? What happens if the spins take value in a Cantor set, or in the $p$-adic numbers? 
The same questions can be asked in $\mathbb{Z}^2$.

%\noindent
%{\em Question:} on the $\Z^2$ grid (or on a planar hyperbolic  lattice) assume the spins at distance $n$ from the origin can get one of $f(n)$, $k 2\pi/ f(n)$ roots of unity  in %$\mathbb{S}^1$. What is the  the critical growth of $f(n)$, for existence of  a magnetization phase?

\end{document}